\documentclass[12pt]{amsart}
\usepackage{amsmath,amscd,amssymb,latexsym}
\usepackage{graphicx}
\input{epsf}
\newtheorem{theorem}{Theorem}
\newtheorem{lemma}{Lemma}
\newtheorem{cor}{Corollary}
\newtheorem{prop}{Proposition}
\def\a{{\alpha}}

\def\d{{\delta}}

\def\s{{\sigma}}

\begin{document}

\title
{The Kauffman Polynomials of 2-bridge Knots}
\author{Bin Lu}
\address{Department of Mathematics \& Statistics\\California State University Sacramento\\6000 J Street\\ Sacramento, CA 95819}
\email{binlu@csus.edu}

\author{Jianyuan K. Zhong}
\address{Department of Mathematics \& Statistics\\California State University Sacramento\\6000 J Street\\ Sacramento, CA 95819}
\email{kzhong@csus.edu}

\date{\today}
\begin{abstract}
The $2$-bridge knots are a family of knots with bridge number $2$ \cite{A94}
\cite{KA}. In this paper, we compute the Kauffman polynomials of $2$-bridge knots using the Kauffman skein theory and linear
algebra techniques. Our calculation can be easily carried out
using Mathematica, Maple, Mathcad, etc. 
\end{abstract}

\maketitle
\section{Introduction}


The $2$-bridge knots (or links) are a family of knots with bridge
number $2$. A $2$-bridge knot (link) has at most $2$ components. Except for the knot $8_5$, the first 25 knots in the Rolfsen Knot Table are $2$-bridge knots. A $2$-bridge knot is also called a rational knot because it can be obtained as the numerator or denominator closure of a rational tangle. The rich mathematical aspects of $2$-bridge knots can be found in many references such as \cite{B}, \cite{GK}, \cite{K}, \cite{KA}, \cite{KL}, \cite{R} and \cite{PY}.
The regular diagram $D$ of a $2$-bridge knot can be drawn as
follows \cite{KA}.
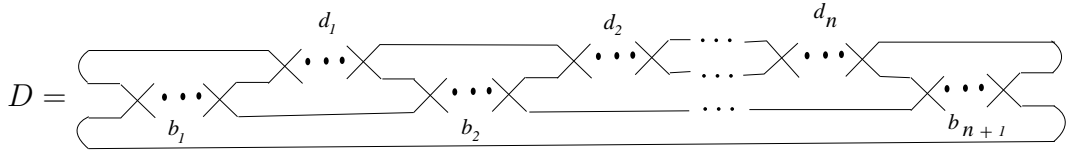
$$D=\raisebox{-6mm}{\epsfysize0.8in\epsffile{2knot.ai}}$$
In the diagram, $d_1, d_2,\cdots, d_n, b_1, b_2,\cdots, b_{n+1}$ are nonzero integers whose absolute values indicate the number of crossings.  
By an isotopy, the diagram $D$ can also be drawn as
$$D=\raisebox{-20mm}{\epsfysize2in\epsffile{2knot1.ai}}$$
The continued fraction notation for $D$ is $F(D):=[ b_1,
d_1, b_2, d_2,\cdots,d_n, b_{n+1}]$ \cite{KL}. We will work with the diagram above to
calculate the Kauffman polynomials of $2$-bridge knots. Specific information about twists (crossings) are necessary to
identify a $2$-bridge knot. Locally, there are four possibilities:
$$\raisebox{-0mm}{\epsfysize0in\epsffile{ltwist.ai}}$$
\hspace{1.5cm}(a) Horizontal left-hand twist \hspace{1.5cm}   (b) Vertical left-hand twist
$$\raisebox{-0mm}{\epsfysize0in\epsffile{rtwist.ai}}$$
\hspace{1.5cm}(c) Horizontal right-hand twist   \hspace{1.3cm}   (d) Vertical right-hand twist

\vspace{0.2cm}

We choose left-hand twists to be positive twists, and right-hand
twists to be negative twists. We will use a positive integer in
the regular diagram to indicate the number of crossings in a
left-hand twist, and a negative integer
to indicate the number of crossings in a right-hand twist.
For example, the Whitehead link is a
2-component 2-bridge link with a diagram given by
$$W=\quad \raisebox{-11mm}{\epsfysize0in
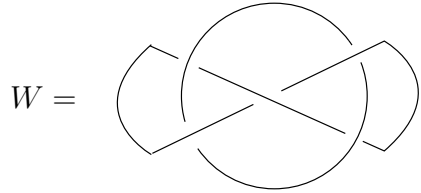
\epsffile{whitehead.ai}}$$
by an isotopy, it can be drawn as
$$W=\quad \raisebox{-17mm}{\epsfysize1.5in
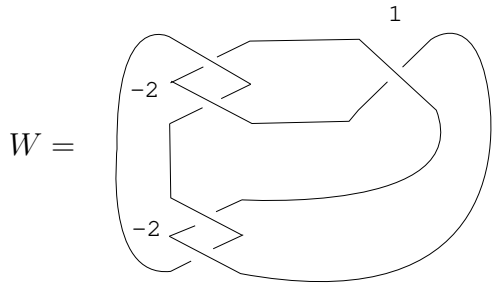
\epsffile{whitehead2.ai}}$$
The continued fraction notation for the Whitehead link is $F(W)=[-2,1,-2]$.

Let ${\mathbb Q}(\alpha, s)$ be the field of rational functions in
$\alpha, s$. By a {\it  framed link} we mean an unoriented link
equipped with a nonsingular normal vector field up to homotopy. The
links described by figures in this paper will be assigned the
vertical framing pointing towards the reader.

There are various versions of {\it the Kauffman polynomial} in the
literature \cite{MT}. Here the Kauffman polynomial of a knot or link is the
unique two-variable rational function in $\alpha, s$ that satisfies the
following Kauffman skein relations:
$$
{\text{(i)}}\quad \raisebox{-3mm}{\epsfxsize.3in\epsffile{left1.ai}}\quad
        -\quad\raisebox{-3mm}{\epsfxsize.3in\epsffile{right1.ai}}\quad
=\quad (\ s - \
s^{-1})\left(\quad\raisebox{-3mm}{\epsfxsize.3in\epsffile{parra1.ai}}-\quad\raisebox{-3mm}{\epsfxsize.3in\epsffile{parra2.ai}}\quad\right
)\quad ,
$$
$$
{\text{(ii)}}\quad \raisebox{-3mm}{\epsfxsize.35in\epsffile{framel1.ai}}
        \quad=\quad \alpha \quad \raisebox{-3mm}{\epsfysize.3in\epsffile{orline1.ai}}\quad  ,
$$
$$
{\text{(iii)}}\quad L\ \sqcup
\raisebox{-2mm}{\epsfysize.3in\epsffile{unknot1.ai}}\quad =\quad
\delta\quad L\quad,
$$
where $\delta=\left({\dfrac{\alpha-\alpha^{-1}} {\ s - \
s^{-1}}}+1\right)$. We denote the Kauffman polynomial of a link
$L$ by $<L>$. Relation (iii) follows from the first two when $L$
is nonempty. Relations (i) and (ii) are local relations,
 except where shown, diagrams in each are identical. A trivial
closed curve in (iii) is a curve which contains no crossing and is null-homotopic.
We normalize the Kauffman polynomial of
the empty link $\emptyset$ to be $1$.

In section 2, we study the Kauffman skein space of the 3-ball $B^3$ with possible boundary points. In section 3, we define linear skein maps on the Kauffman skein space of the 3-ball $B^3$ with four boundary points and compute the matrices of these linear maps. In section 4, we present our main theorem of calculating the Kauffman polynomial of a 2-bridge knot by decomposing it as compositions of linear skein maps from section 3.
In section 5, we calculate the Kauffman polynomial of the Whitehead link as an example.

\section{The Kauffman skein space of the $3$-ball $B^3$}

\subsection{The Kauffman skein space of the $3$-ball $B^3$}
The Kauffman skein space \cite{BB2000} \cite{ZLu} of the $3$-ball
$B^3$, denoted by $K(B^3)$, is the ${\mathbb Q}(\alpha, s)$-space
freely generated by framed isotopic links $L$ in $B^3$ quotient by
the subspace generated by the Kauffman skein relations. Given any
link $L$ in $B^3$, it can be simplified to $<L>\emptyset$ by
applying the Kauffman skein relation, where $<L>$ is the Kauffman
polynomial of $L$. Hence the Kauffman skein space $K(B^3)$ is
generated by the empty link $\emptyset$.

\subsection{The Kauffman skein space of $B^3$ with four boundary points}

We place a distinguished set of four coplanar points \{$N, E,
S,W$\} on the sphere $S^2$, the boundary of the $3$-ball $B^3$. A
link in $(B^3, NESW)$ is a collection of closed curves and arcs
joining the distinguished boundary points $N, E, S,W$. Two links
are equivalent if one can be obtained from the other by 
isotopy. We define the Kauffman skein space $K(B^3, NESW)$ to be
the ${\mathbb Q}(\alpha, s)$-space freely generated by framed
links $L$ in $(B^3,S^2)$ such that $L \cap S^2=\partial L=\{$N, E,
S,W$\}$, considered up to an ambient isotopy fixing $S^2$, quotient
by the subspace generated by the Kauffman skein relations. A
 skein element in $K(B^3, NESW)$ is illustrated below.
$$\raisebox{-8mm}{\epsfysize0.8in\epsffile{diagramb3.ai}}$$
There are two natural multilinear multiplication operations  in $K(B^3, NESW)$:
\begin{enumerate}

\item {\bf Concatenation.} By stacking the first on top of the second through gluing points $W,S$ in the first with
$N,E$ in the second,
$$\odot: \quad \raisebox{-8mm}{\epsfysize0.8in\epsffile{diagramb3.ai}}\odot \raisebox{-8mm}{\epsfysize0.8in\epsffile{diagramb33.ai}}\quad =\quad \raisebox{-8mm}{\epsfysize0.8in\epsffile{diagramb32.ai}}\quad;$$
\item {\bf Juxtaposition.} By putting two skein elements next to each other through gluing points $E,S$ in the first  with $N,W$ in
the second,
$$\otimes: \quad \raisebox{-8mm}{\epsfysize0.8in\epsffile{diagramb3.ai}}\otimes \raisebox{-8mm}{\epsfysize0.8in\epsffile{diagramb33.ai}}\quad =\quad \raisebox{-8mm}{\epsfysize0.8in\epsffile{diagramb31.ai}}\quad.$$
\end{enumerate}

Note that the skein element $\raisebox{-3mm}{\epsfxsize.3in\epsffile{parra1.ai}}$ is the identity with respect to the $\odot$ operation, and the skein element $\raisebox{-3mm}{\epsfxsize.3in\epsffile{parra2.ai}}$ is the identity with respect to the $\otimes$ operation.

The Kauffman skein space $K(B^3, NESW)$ is 3-dimensional and has a basis
\{${e_1}$, ${e_2}$, $e_3$\} \cite{BB2000} given by

$\hspace{1cm}{e_1}=\displaystyle{\frac{1}{s+s^{-1}}}\left(s^{-1}\quad \raisebox{-3mm}{\epsfxsize.3in\epsffile{parra1.ai}}\quad + \quad \raisebox{-3mm}{\epsfxsize.3in\epsffile{left1.ai}}\quad-(\delta^{-1}s^{-1}+ \delta^{-1}\alpha^{-1}) \quad \raisebox{-3mm}{\epsfxsize.3in\epsffile{parra2.ai}}\right);$
\vspace{0.2cm}

$\hspace{1cm}{e_2}=\displaystyle{\frac{1}{s+s^{-1}}}\left(s\quad \raisebox{-3mm}{\epsfxsize.3in\epsffile{parra1.ai}}\quad - \quad \raisebox{-3mm}{\epsfxsize.3in\epsffile{left1.ai}}\quad+ (-\delta^{-1}s+ \delta^{-1}\alpha^{-1}) \quad \raisebox{-3mm}{\epsfxsize.3in\epsffile{parra2.ai}}\right);$
\vspace{0.2cm}

$\hspace{1cm}{e_3}=\d^{-1}\quad\raisebox{-3mm}{\epsfxsize0in\epsffile{parra2.ai}}.$

In the remaining part of this section, we study properties of
these basis elements which are crucial in
constructing our calculation techniques.
\begin{prop} With respect to the $\odot$ operation,
\begin{enumerate}
\item the basis elements ${e_1}$, ${e_2}$, $e_3$ are
orthogonal,  i.e.,
$${e_1}\odot {e_2}={e_2}\odot {e_1}=0,\quad {e_1}\odot {e_3}={e_3}\odot{e_1}=0,\quad {e_2}\odot {e_3}={e_3}\odot {e_2}=0;$$
\item the basis elements ${e_1}$,  ${e_2}$ and $e_3$ are
idempotents, i.e., ${e_1}\odot{e_1}={e_1}$,
${e_2}\odot{e_2}={e_2}$, $e_3\odot e_3= e_3$; \item the basis
elements ${e_1}$, ${e_2}$, $e_3$ add to the identity with respect to the $\odot$ operation, i.e.,
${e_1}+ {e_2}
+{e_3}=\raisebox{-3mm}{\epsfxsize.3in\epsffile{parra1.ai}}$;
\item let $\s=\raisebox{-3mm}{\epsfxsize.3in\epsffile{left1.ai}}\quad$, then 
$\s \odot {e_1}={e_1}\odot \s=s{e_1}, $ $\s \odot {e_2}={e_2}\odot
\s=-s^{-1}{e_2},$ $\s \odot {e_3}={e_3}\odot \s=\a^{-1}{e_3}.$ It follows that 
$\s^{-1} \odot {e_1}={e_1}\odot \s^{-1}=s^{-1}{e_1}, $ $\s^{-1} \odot {e_2}={e_2}\odot
\s^{-1}=-s{e_2},$ $\s^{-1} \odot {e_3}={e_3}\odot \s^{-1}=\a {e_3}.$
\end{enumerate}
\end{prop}
Let $\s_\odot^n$ represent $n$ copies of $\s$ multiplied through
the ``$\odot$" multiplication structure, it follows that
$\s_\odot^n \odot {e_1}={e_1}\odot \s_\odot^n=s^n{e_1}$,
$\s_\odot^n\odot {e_2}={e_2}\odot \s_\odot^n=(-s^{-1})^{n}{e_2},$
$\s_\odot^n \odot {e_3}={e_3}\odot \s_\odot^n=\a^{-n}{e_3}$.
\begin{proof}
The proofs follow by the linearity of the $\odot$ operation and (repeatedly) applying the Kauffman skein
relations and substituting
$\delta=\left({\dfrac{\alpha-\alpha^{-1}} {\ s - \
s^{-1}}}+1\right)$. Here we show ${e_1}\odot {e_2}=0$ as an example.

${e_1}\odot {e_2}=\displaystyle{\frac{1}{(s+s^{-1})^2}}
\Big(s^{-1}\raisebox{-3mm}{\epsfxsize.3in\epsffile{parra1.ai}} +  \raisebox{-3mm}{\epsfxsize.3in\epsffile{left1.ai}}-(\delta^{-1}s^{-1}+ \delta^{-1}\alpha^{-1}) \raisebox{-3mm}{\epsfxsize.3in\epsffile{parra2.ai}}\Big)\odot \Big(s \raisebox{-3mm}{\epsfxsize.3in\epsffile{parra1.ai}}- \raisebox{-3mm}{\epsfxsize.3in\epsffile{left1.ai}}+ (-\delta^{-1}s+ \delta^{-1}\alpha^{-1})\raisebox{-3mm}{\epsfxsize.3in\epsffile{parra2.ai}}\Big)$

$=\displaystyle{\frac{1}{(s+s^{-1})^2}}\Big(\raisebox{-3mm}{\epsfxsize.3in\epsffile{parra1.ai}} -s^{-1}  \raisebox{-3mm}{\epsfxsize.3in\epsffile{left1.ai}}+ s^{-1}(-\delta^{-1}s+ \delta^{-1}\alpha^{-1}) \raisebox{-3mm}{\epsfxsize.3in\epsffile{parra2.ai}}+ s \raisebox{-3mm}{\epsfxsize.3in\epsffile{left1.ai}} - \raisebox{-3mm}{\epsfxsize.3in\epsffile{left1.ai}}\odot \raisebox{-3mm}{\epsfxsize.3in\epsffile{left1.ai}}+ (-\delta^{-1}s+ \delta^{-1}\alpha^{-1}) \alpha^{-1}\raisebox{-3mm}{\epsfxsize.3in\epsffile{parra2.ai}}-(\delta^{-1}s^{-1}+ \delta^{-1}\alpha^{-1})s  \raisebox{-3mm}{\epsfxsize.3in\epsffile{parra2.ai}}+ (\delta^{-1}s^{-1}+ \delta^{-1}\alpha^{-1})\a^{-1} \raisebox{-3mm}{\epsfxsize.3in\epsffile{parra2.ai}} - (\delta^{-1}s^{-1}+ \delta^{-1}\alpha^{-1})(-\delta^{-1}s+ \delta^{-1}\alpha^{-1})\d \quad \raisebox{-3mm}{\epsfxsize.3in\epsffile{parra2.ai}}\Big)=0$.
\end{proof}

If we rotate the basis elements $e_1, e_2, e_3$ in the plane by $90^\circ$, we obtain another basis for $K(B^3, NESW)$.  We present the basis elements ${e_1}$, ${e_2}$, $e_3$
using subscripts $h$ (vs $v$) to indicate the basis elements after (vs before) the rotation:

${e_1}_v=e_1=\displaystyle{\frac{1}{s+s^{-1}}}\left(s^{-1}\quad \raisebox{-3mm}{\epsfxsize.3in\epsffile{parra1.ai}} + \quad \raisebox{-3mm}{\epsfxsize.3in\epsffile{left1.ai}}-(\delta^{-1}s^{-1}+ \delta^{-1}\alpha^{-1}) \quad \raisebox{-3mm}{\epsfxsize.3in\epsffile{parra2.ai}}\right);$
\vspace{0.2cm}

${e_2}_v=e_2=\displaystyle{\frac{1}{s+s^{-1}}}\left(s\quad \raisebox{-3mm}{\epsfxsize.3in\epsffile{parra1.ai}} - \quad \raisebox{-3mm}{\epsfxsize.3in\epsffile{left1.ai}}\quad+(-\delta^{-1}s+ \delta^{-1}\alpha^{-1}) \quad \raisebox{-3mm}{\epsfxsize.3in\epsffile{parra2.ai}}\right);$
\vspace{0.2cm}

${e_3}_v=e_3=\delta^{-1}\quad\raisebox{-3mm}{\epsfxsize0in\epsffile{parra2.ai}}.$
\vspace{0.2cm}

${e_1}_h=\displaystyle{\frac{1}{s+s^{-1}}}\left(s^{-1}\quad
\raisebox{-3mm}{\epsfxsize.3in\epsffile{parra2.ai}} + \quad
\raisebox{-3mm}{\epsfxsize.3in\epsffile{right1.ai}}-
(\delta^{-1}s^{-1}+ \delta^{-1}\alpha^{-1}) \quad
\raisebox{-3mm}{\epsfxsize.3in\epsffile{parra1.ai}}\right);$
\vspace{0.2cm}

${e_2}_h=\displaystyle{\frac{1}{s+s^{-1}}}\left(s\quad
\raisebox{-3mm}{\epsfxsize.3in\epsffile{parra2.ai}}\quad -
\raisebox{-3mm}{\epsfxsize.3in\epsffile{right1.ai}}\quad+(-\delta^{-1}s+
\delta^{-1}\alpha^{-1}) \quad
\raisebox{-3mm}{\epsfxsize.3in\epsffile{parra1.ai}}\right);$
\vspace{0.2cm}

${e_3}_h=\delta^{-1}\quad\raisebox{-3mm}{\epsfxsize0in\epsffile{parra1.ai}}.$

With respect to the $\otimes$ operation and the basis elements ${e_1}_h, {e_2}_h, {e_3}_h $, similar properties of the basis
elements given in Proposition 1 still hold, which we state as a corollary below.
\begin{cor}
\begin{enumerate}
\item ${e_1}_h\otimes {e_2}_h={e_2}_h\otimes {e_1}_h=0,\quad
{e_1}_h\otimes {e_3}_h={e_3}_h\otimes{e_1}_h=0,$

$\quad {e_2}_h\otimes {e_3}_h={e_3}_h\otimes {e_2}_h=0;$
\item
${e_1}_h\otimes{e_1}_h={e_1}_h$,
${e_2}_h\otimes{e_2}_h={e_2}_h$, ${e_3}_h\otimes {e_3}_h=
{e_3}_h$; \item
$\raisebox{-3mm}{\epsfxsize.3in\epsffile{parra2.ai}}={e_1}_h+
{e_2}_h +{e_3}_h$;
\item let $\s_h=\raisebox{-3mm}{\epsfxsize.3in\epsffile{right1.ai}}$, then 
$\s_h \otimes {e_1}_h={e_1}_h\otimes \s_h=s{e_1}_h, $
${\s_h}\otimes {e_2}_h={e_2}_h\otimes \s_h=-s^{-1}{e_2}_h, $
${\s_h} \otimes {e_3}_h={e_3}_h\otimes
{\s_h}=\a^{-1}{e_3}_h.$
\end{enumerate}
\end{cor}
It follows that ${\s_h}_\otimes^n \otimes {e_1}_h={e_1}_h\otimes
{\s_h}_\otimes^n=s^n{e_1}_h$, ${\s_h}_\otimes^n\otimes
{e_2}_h={e_2}_h\otimes {\s_h}_\otimes^n=(-s^{-1})^{n}{e_2}_h,$
${\s_h}_\otimes^n \otimes {e_3}_h={e_3}_h\otimes
{\s_h}_\otimes^n=\a^{-n}{e_3}_h$, where ${\s_h}_\otimes^n$
represents $n$ copies of $\s_h$ multiplied through the
``$\otimes$" operation.

The following are additional properties of the basis elements ${e_1}_v, {e_2}_v, {e_3}_v $, ${e_1}_h, {e_2}_h, {e_3}_h $ with respect to the $\otimes$ operation.
\begin{prop}
\begin{enumerate}
\item $\displaystyle{{e_1}_h \otimes
{e_1}_v={e_1}_v \otimes
{e_1}_h=\frac{1}{s+s^{-1}}(s^{-1}-\delta^{-1}s^{-1}-
\delta^{-1}\alpha^{-1}){e_1}_h}$;

\item $\displaystyle{{e_1}_h \otimes
{e_2}_v={e_2}_v \otimes
{e_1}_h=\frac{1}{s+s^{-1}}(-s^{-1}-\delta^{-1}s+
\delta^{-1}\alpha^{-1}){e_1}_h}$;

\item $\displaystyle{{e_1}_h \otimes
{e_3}_v={e_3}_v \otimes
{e_1}_h=\delta^{-1}{e_1}_h}$;

\item $\displaystyle{{e_2}_h \otimes
{e_1}_v={e_1}_v \otimes
{e_2}_h=\frac{1}{s+s^{-1}}(-s-\delta^{-1}s^{-1}-
\delta^{-1}\alpha^{-1}){e_2}_h}$;

\item $\displaystyle{{e_2}_h \otimes
{e_2}_v={e_2}_v \otimes
{e_2}_h=\frac{1}{s+s^{-1}}(s-\delta^{-1}s+
\delta^{-1}\alpha^{-1}){e_2}_h}$;

\item $\displaystyle{{e_2}_h \otimes
{e_3}_v={e_3}_v \otimes
{e_2}_h=\delta^{-1}{e_2}_h}$;

\item $\displaystyle{{e_3}_h \otimes
{e_1}_v={e_1}_v \otimes
{e_3}_h=\frac{1}{s+s^{-1}}(s^{-1}\d+\a-\delta^{-1}s^{-1}-
\delta^{-1}\alpha^{-1}){e_3}_h}$;

\item $\displaystyle{{e_3}_h \otimes
{e_2}_v={e_2}_v \otimes
{e_3}_h=\frac{1}{s+s^{-1}}(s\d -\a-\delta^{-1}s+
\delta^{-1}\alpha^{-1}){e_3}_h}$;

\item $\displaystyle{{e_3}_h \otimes
{e_3}_v={e_3}_v \otimes
{e_3}_h=\delta^{-1}{e_3}_h}$;

\end{enumerate}
\end{prop}


\begin{proof} Here we prove (1) as an example, (2)-(9) can be proved in a similar fashion.

$\displaystyle{{e_1}_h \otimes {e_1}_v}= \displaystyle{\frac{1}{s+s^{-1}}{e_1}_h \otimes \left(s^{-1}\quad \raisebox{-3mm}{\epsfxsize.3in\epsffile{parra1.ai}}\quad + \quad \raisebox{-3mm}{\epsfxsize.3in\epsffile{left1.ai}}\quad-(\delta^{-1}s^{-1}+ \delta^{-1}\alpha^{-1}) \quad \raisebox{-3mm}{\epsfxsize.3in\epsffile{parra2.ai}}\right)}$

$=\displaystyle{\frac{1}{s+s^{-1}}\left(s^{-1}{e_1}_h \otimes \quad \raisebox{-3mm}{\epsfxsize.3in\epsffile{parra1.ai}}\quad +{e_1}_h \otimes {\s_h}^{-1}-(\delta^{-1}s^{-1}+ \delta^{-1}\alpha^{-1}){e_1}_h \otimes \raisebox{-3mm}{\epsfxsize.3in\epsffile{parra2.ai}}\right)}$

$=\displaystyle{\frac{1}{s+s^{-1}}\left(0+s^{-1}{e_1}_h-(\delta^{-1}s^{-1}+ \delta^{-1}\alpha^{-1}){e_1}_h\right)}$

$=\displaystyle{\frac{1}{s+s^{-1}}(s^{-1}-\delta^{-1}s^{-1}-
\delta^{-1}\alpha^{-1}){e_1}_h}$.
\end{proof}

Let $M$ be the $3\times 3$ matrix given by
$$M=\left(\begin{matrix}{\displaystyle{\frac{1}{s+s^{-1}}(s^{-1}-\delta^{-1}s^{-1}-
\delta^{-1}\alpha^{-1})}} &
{\displaystyle{\frac{1}{s+s^{-1}}(-s^{-1}-\delta^{-1}s+
\delta^{-1}\alpha^{-1})}} & \d^{-1}\\ & &
\\{\displaystyle{\frac{1}{s+s^{-1}}(-s-\delta^{-1}s^{-1}-
\delta^{-1}\alpha^{-1})}} &
{\displaystyle{\frac{1}{s+s^{-1}}(s-\delta^{-1}s+
\delta^{-1}\alpha^{-1})}} & \d^{-1}\\ & & \\
{\displaystyle{\frac{1}{s+s^{-1}}(s^{-1}\d+\a-\d^{-1}s^{-1}-
\delta^{-1}\alpha^{-1})}} &
{\displaystyle{\frac{1}{s+s^{-1}}(s\d-\a-\d^{-1}s+
\delta^{-1}\alpha^{-1})}} & \d^{-1} \end{matrix}\right )$$
\noindent {\bf Remark.} The entries of $M=(m_{ij})$ are the
coefficients in Proposition 2 (1)-(9), where ${e_i}_h\otimes
{e_j}_v=m_{ij}{e_i}_h$ for $1\leq i,j\leq 3$. Notice that the matrix
$M$ is the base change matrix between the basis
$\{{e_1}_h, {e_2}_h, {e_3}_h\}$ and $\{{e_1}_v, {e_2}_v,
{e_3}_v\}$, i.e.,
$$({e_1}_v,{e_2}_v,{e_3}_v)=({e_1}_h, {e_2}_h, {e_3}_h)M, \quad ({e_1}_h,{e_2}_h,{e_3}_h)=({e_1}_v, {e_2}_v, {e_3}_v)M.$$
It follows that $M^2=I$, the $3\times 3$ identity matrix.

\noindent {\bf Remark: } If we change $\otimes$ to $\odot$ and exchange the subscripts $v$ and $h$ in
identities (1)-(9) in Proposition 2, the identities  still hold. We state these in the next corollary.
\begin{cor}
\begin{enumerate}
\item $\displaystyle{{e_1}_v \odot
{e_1}_h={e_1}_h \odot
{e_1}_v=\frac{1}{s+s^{-1}}(s^{-1}-\delta^{-1}s^{-1}-
\delta^{-1}\alpha^{-1}){e_1}_v}$;

\item $\displaystyle{{e_1}_v \odot
{e_2}_h={e_2}_h \odot
{e_1}_v=\frac{1}{s+s^{-1}}(-s^{-1}-\delta^{-1}s+
\delta^{-1}\alpha^{-1}){e_1}_v}$;

\item $\displaystyle{{e_1}_v \odot
{e_3}_h={e_3}_h \odot
{e_1}_v=\delta^{-1}{e_1}_v}$;

\item $\displaystyle{{e_2}_v \odot
{e_1}_h={e_1}_h \odot
{e_2}_v=\frac{1}{s+s^{-1}}(-s-\delta^{-1}s^{-1}-
\delta^{-1}\alpha^{-1}){e_2}_v}$;

\item $\displaystyle{{e_2}_v \odot
{e_2}_h={e_2}_h \odot
{e_2}_v=\frac{1}{s+s^{-1}}(s-\delta^{-1}s+
\delta^{-1}\alpha^{-1}){e_2}_v}$;

\item $\displaystyle{{e_2}_v \odot
{e_3}_h={e_3}_h \odot
{e_2}_v=\delta^{-1}{e_2}_v}$;

\item $\displaystyle{{e_3}_v \odot
{e_1}_h={e_1}_h \odot
{e_3}_v=\frac{1}{s+s^{-1}}(s^{-1}\d+\a-\delta^{-1}s^{-1}-
\delta^{-1}\alpha^{-1}){e_3}_v}$;

\item $\displaystyle{{e_3}_v \odot
{e_2}_h={e_2}_h \odot
{e_3}_v=\frac{1}{s+s^{-1}}(s\d -\a-\delta^{-1}s+
\delta^{-1}\alpha^{-1}){e_3}_v}$;

\item $\displaystyle{{e_3}_v \odot
{e_3}_h={e_3}_h \odot
{e_3}_v=\delta^{-1}{e_3}_v}$.

\end{enumerate}
\end{cor}


\section{Linear Skein Maps on $K(B^3, NESW)$ and their matrices}

A wiring of a space $F$ into another space $F'$ is a choice of
inclusion of $F$ into $F'$ and a choice of a set of fixed curves
and arcs in $F'-F$. The wiring of $F$ into $F'$ induces a well-defined linear map from the skein space $K(F)$ to $K(F')$ \cite{HM}. In this section we'll consider four wirings of
$B^3$ into itself, three of which induce linear skein maps $K(B^3,
NESW)\to K(B^3, NESW)$, while the fourth one induces a linear skein map
$K(B^3, NESW)\to K(B^3)$. Since $K(B^3, NESW)$ and $K(B^3)$ are
vector spaces over ${\mathbb Q}(\alpha, s)$, these linear
maps are linear transformations of vector spaces.
In the following, we choose $\{{e_1}_h, {e_2}_h, {e_3}_h\}$ as the basis of $K(B^3, NESW)$ 
and represent these linear transformations by
matrices with respect to this basis.

\subsection{The linear map $B_1(b_1)$ and the matrix $B_1(b_1)$}

Let $b_1$ be a nonzero integer, the linear map $B_1(b_1): K(B^3, NESW)\to
K(B^3, NESW),$ is induced by the following wiring, also called $B_1(b_1)$, for convenience 
$$B_1(b_1):\quad \raisebox{-8mm}{\epsfysize0.8in
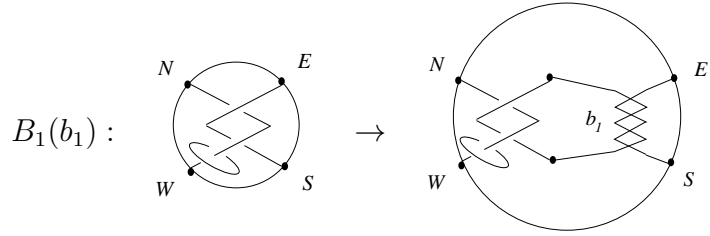
\epsffile{diagramb3.ai}}\quad \to \quad \raisebox{-12mm}{\epsfysize1.2in\epsffile{b1map.ai}}$$
where $b_1$ indicates the number of crossings, it is positive if the crossings form left-hand twists, it is negative if the crossings form right-hand twists.

\begin{lemma}$B_1(b_1) (x{e_1}_h+y{e_2}_h+z{e_3}_h)=
x(m_{11}s^{b_1}+m_{12}(-s^{-1})^{b_1}+m_{13}\a^{-b_1}){e_1}_h+y(m_{21}s^{b_1}+m_{22}(-s^{-1})^{b_1}+m_{23}\a^{-b_1}){e_2}_h+z(m_{31}s^{b_1}+m_{32}(-s^{-1})^{b_1}+m_{33}\a^{-b_1}){e_3}_h$.

\end{lemma}
\begin{proof}
Let $s\in K(B^3, NESW)\to
K(B^3, NESW)$, then $B_1(b_1) (s)=s\otimes (\s_\odot^{b_1})$. Note that $\s_\odot^{b_1}=\s_\odot^{b_1}\odot \left(\quad\raisebox{-3mm}{\epsfxsize.3in\epsffile{parra1.ai}}\right)=\s_\odot^{b_1}\odot \left({e_1}_v +{e_2}_v+{e_3}_v\right)=\s_\odot^{b_1}\odot {e_1}_v +\s_\odot^{b_1}\odot {e_2}_v+\s_\odot^{b_1}\odot {e_3}_v=s^{b_1}{e_1}_v +(-s^{-1})^{b_1}{e_2}_v+\a^{-b_1}{e_3}_v$ from Proposition 1.

Now $B_1(b_1) (x{e_1}_h+y{e_2}_h+z{e_3}_h)=(x{e_1}_h+y{e_2}_h+z{e_3}_h)\otimes \s_\odot^{b_1}$

$=(x{e_1}_h+y{e_2}_h+z{e_3}_h)\otimes \left(s^{b_1}{e_1}_v +(-s^{-1})^{b_1}{e_2}_v+\a^{-b_1}{e_3}_v\right)$

$=x{e_1}_h\otimes \left(s^{b_1}{e_1}_v +(-s^{-1})^{b_1}{e_2}_v+\a^{-b_1}{e_3}_v\right) +y{e_2}_h\otimes \left(s^{b_1}{e_1}_v +(-s^{-1})^{b_1}{e_2}_v+\a^{-b_1}{e_3}_v\right) +z{e_3}_h\otimes \left(s^{b_1}{e_1}_v +(-s^{-1})^{b_1}{e_2}_v+\a^{-b_1}{e_3}_v\right)$

$=xs^{b_1}{e_1}_h\otimes {e_1}_v +x(-s^{-1})^{b_1}{e_1}_h\otimes {e_2}_v+x\a^{-b_1}{e_1}_h\otimes {e_3}_v+ys^{b_1}{e_2}_h\otimes {e_1}_v +y(-s^{-1})^{b_1}{e_2}_h\otimes {e_2}_v+y\a^{-b_1}{e_2}_h\otimes {e_3}_v+zs^{b_1}{e_3}_h\otimes {e_1}_v +z(-s^{-1})^{b_1}{e_3}_h\otimes {e_2}_v+z\a^{-b_1}{e_3}_h\otimes {e_3}_v$

$=xs^{b_1}m_{11}{e_1}_h +x(-s^{-1})^{b_1}m_{12}{e_1}_h+x\a^{-b_1}m_{13}{e_1}_h+ys^{b_1}m_{21}{e_2}_h +y(-s^{-1})^{b_1}m_{22}{e_2}_h+y\a^{-b_1}m_{23}{e_2}_h+zs^{b_1}m_{31}{e_3}_h+z(-s^{-1})^{b_1}m_{32}{e_3}_h+z\a^{-b_1}m_{33}{e_3}_h$

$=x(m_{11}s^{b_1}+m_{12}(-s^{-1})^{b_1}+m_{13}){e_1}_h+y(m_{21}s^{b_1}+m_{22}(-s^{-1})^{b_1}+m_{23}\a^{-b_1}){e_2}_h+z(m_{31}s^{b_1}+m_{32}(-s^{-1})^{b_1}+m_{33}\a^{-b_1}){e_3}_h.$
\end{proof}

We define the corresponding matrix $B_1(b_1)$ as
$B_1(b_1)=(b_{ij}), 1\leq i,j\leq 3,$
where 
$$b_{ij}=\left\{\begin{array}{cc}
m_{11}s^{b_1}+m_{12}(-s^{-1})^{b_1}+m_{13}\a^{-b_1}  & \mbox{if $i=j=1$}\\ & \\
m_{21}s^{b_1}+m_{22}(-s^{-1})^{b_1}+m_{23}\a^{-b_1} & \mbox{if $i=j=2$}\\ &\\
m_{31}s^{b_1}+m_{32}(-s^{-1})^{b_1}+m_{33}\a^{-b_1} & \mbox{if $i=j=3$} \\ &\\
0 & \mbox{Otherwise.} \\
\end{array}
\right.$$


\subsection{The linear map $D(d_i)$ and the matrix $D(d_i)$}

Let $d_i$ be a nonzero integer, the linear map $D(d_i): K(B^3, NESW)\to
K(B^3, NESW)$ is induced by the wiring
$$D(d_i):\quad\raisebox{-8mm}{\epsfysize0.8in
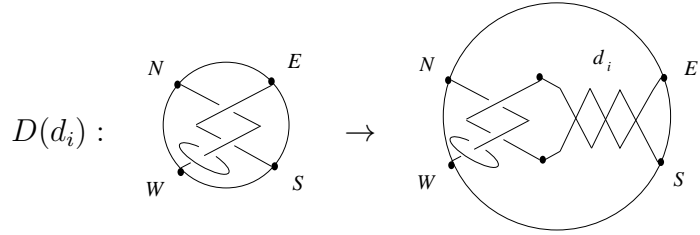
\epsffile{diagramb3.ai}}\quad \to \quad \raisebox{-12mm}{\epsfysize1.2in\epsffile{dimap.ai}}$$
Similarly $d_i$ indicates the number of crossings, it is positive if the crossings form left-hand twists, it is negative if the crossings form right-hand twists.

\begin{lemma}
$D(d_i) (x{e_1}_h+y{e_2}_h+z{e_3}_h)=xs^{d_i}{e_1}_h
+y(-s^{-1})^{d_i}{e_2}_h + z\a^{-d_i}{e_3}_h$
$$=({e_1}_h\ {e_2}_h\ {e_3}_h)\left(\begin{matrix}s^{d_i} & 0 & 0\\ 0 & (-s^{-1})^{d_i} & 0\\ 0 & 0 & \a^{-d_i} \end{matrix}\right)\left(\begin{matrix}x\\ y\\ z \end{matrix}\right ).$$
\end{lemma}
\begin{proof} Note that
${\s_h}_\otimes^{d_i}=
{\s_h}_\otimes^{d_i} \otimes
\left(\quad\raisebox{-3mm}{\epsfxsize.3in\epsffile{parra2.ai}}\right)
={\s_h}_\otimes^{d_i} \otimes \left({e_1}_h
+{e_2}_h+{e_3}_h\right)$

$=s^{d_i}{e_1}_h
+(-s^{-1})^{d_i}{e_2}_h + \a^{-d_i}{e_3}_h$ by the idempotent properties of the basis elements.

Now by substitution, $D(d_i) (x{e_1}_h+y{e_2}_h+z{e_3}_h)=(x{e_1}_h+y{e_2}_h+z{e_3}_h) \otimes ({\s_h}_\otimes^{d_i})$

$= (x{e_1}_h+y{e_2}_h+z{e_3}_h) \otimes (s^{d_i}{e_1}_h
+(-s^{-1})^{d_i}{e_2}_h + \a^{-d_i}{e_3}_h)$

$=xs^{d_i}{e_1}_h
+y(-s^{-1})^{d_i}{e_2}_h + z\a^{-d_i}{e_3}_h$.

\end{proof}
We define the corresponding matrix $D(d_i)$ as
$$D(d_i)=\left(\begin{matrix}s^{d_i} & 0 & 0\\ 0 & (-s^{-1})^{d_i} & 0\\ 0 & 0 & \a^{-d_i} \end{matrix}\right ).$$

\subsection{The linear map $B(b_i)$ and the matrix $B(b_i)$ }

Let $b_i$ be a nonzero integer, the linear map $B(b_i): K(B^3,
NESW)\to K(B^3, NESW)$ is induced by the wiring

$$B(b_i):\quad\raisebox{-8mm}{\epsfysize0.8in
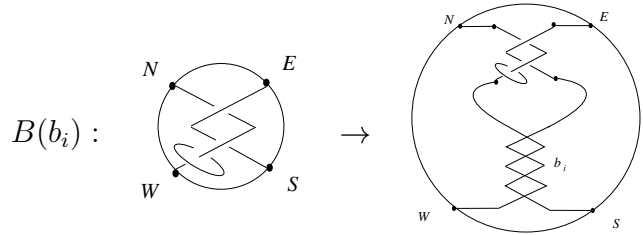
\epsffile{diagramb3.ai}}\quad \to \quad \raisebox{-12mm}{\epsfysize1.2in\epsffile{bimap.ai}}$$
where $b_i$ indicates the number of crossings, it is positive if
the crossings form left-hand twists, it is negative if the
crossings form right-hand twists.

\begin{lemma}
$B(b_i) (x{e_1}_h+y{e_2}_h+z{e_3}_h)=({e_1}_h\
{e_2}_h\ {e_3}_h)M\left(\begin{matrix}s^{b_i} & 0 & 0\\ 0 &
(-s^{-1})^{b_i} & 0\\ 0 & 0 & \a^{-b_i} \end{matrix}\right
)M\left(\begin{matrix}x\\ y\\ z \end{matrix}\right );$
\noindent where $M$ is the base change matrix between the basis
$\{{e_1}_h, {e_2}_h, {e_3}_h\}$ and $\{{e_1}_v, {e_2}_v,
{e_3}_v\}$.
\end{lemma}

\begin{proof}

Note $\s_\odot^{b_i}=\s_\odot^{b_i}\odot \quad\raisebox{-3mm}{\epsfxsize.3in\epsffile{parra1.ai}}=\s_\odot^{b_i}\odot ({e_1}_v +{e_2}_v+{e_3}_v)=s^{b_1}{e_1}_v +(-s^{-1})^{b_1}{e_2}_v+\a^{-b_1}{e_3}_v$,

$B(b_i) (x{e_1}_h+y{e_2}_h+z{e_3}_h)=(x{e_1}_h+y{e_2}_h+z{e_3}_h)\odot (\s_\odot^{b_i})$

$=(x{e_1}_h+y{e_2}_h+z{e_3}_h)\odot \left(s^{b_1}{e_1}_v +(-s^{-1})^{b_1}{e_2}_v+\a^{-b_1}{e_3}_v\right)$

$=\left(x{e_1}_h+y{e_2}_h+z{e_3}_h)\odot \left(s^{b_1}{e_1}_v +(-s^{-1})^{b_1}{e_2}_v+\a^{-b_1}{e_3}_v\right)\right) $

$=(x{e_1}_h+y{e_2}_h+z{e_3}_h)\odot s^{b_1}{e_1}_v + (x{e_1}_h+y{e_2}_h+z{e_3}_h)\odot (-s^{-1})^{b_1}{e_2}_v+(x{e_1}_h+y{e_2}_h+z{e_3}_h)\odot \a^{-b_1}{e_3}_v$

$=(xm_{11}+ym_{12}+zm_{13})s^{b_1}{e_1}_v + (xm_{21}+ym_{22}+zm_{23}) (-s^{-1})^{b_1}{e_2}_v+(xm_{31}+ym_{32}+zm_{33}) \a^{-b_1}{e_3}_v$

$=({e_1}_v\ {e_2}_v\ {e_3}_v)\left(\begin{matrix}s^{b_i} & 0 & 0\\ 0 & (-s^{-1})^{b_i} & 0\\ 0 & 0 & \a^{-b_i} \end{matrix}\right )M\left(\begin{matrix}x\\ y\\ z \end{matrix}\right )$

$=({e_1}_h\ {e_2}_h\ {e_3}_h)M\left(\begin{matrix}s^{b_i} & 0 & 0\\ 0 & (-s^{-1})^{b_i} & 0\\ 0 & 0 & \a^{-b_i} \end{matrix}\right )M\left(\begin{matrix}x\\ y\\ z \end{matrix}\right )$, as $({e_1}_v\ {e_2}_v\ {e_3}_v)=({e_1}_h\ {e_2}_h\ {e_3}_h)M$.
\end{proof}

We define the corresponding matrix $B(b_i)$ as
$$B(b_i)= M\left(\begin{matrix}s^{b_i} & 0 & 0\\ 0 & (-s^{-1})^{b_i} & 0\\ 0 & 0 & \a^{-b_i} \end{matrix}\right )M.$$

\subsection{The closure-map $C$ and the matrix $C$}
Finally the linear  map
$C: K(B^3, NESW)\to K(B^3)$
is induced by the closure wiring:
$$\raisebox{-8mm}{\epsfysize0.8in\epsffile{diagramb3.ai}}\quad \to \quad \raisebox{-12mm}{\epsfysize1.2in\epsffile{cmap.ai}}$$
\begin{lemma}
$$C (x{e_1}_h+y{e_2}_h+z{e_3}_h)=z \d\emptyset,$$
where $\emptyset$ represents the empty link which generates $K(B^3)$.
\end{lemma}
\begin{proof}
The closure of ${e_1}_h$ is zero and the closure of ${e_2}_h$ is also zero by the orthogonal properties. The closure of ${e_3}_h$ can be simplified as $\d^{-1}\d^2 \emptyset=\d\emptyset$.
\end{proof}
We therefore define the matrix $C=(0,0,\d).$

\section{The Kauffman Polynomials of the $2$-bridge knots}

The 2-bridge knot with continuous fraction notation $F(D)=[ b_1,
d_1, b_2, d_2,\cdots,d_n, b_{n+1}]$ is an image of the compositions of wiring maps defined in last section. We summarize our main results in:
\begin{theorem}
Let $F(D)=[ b_1,
d_1, b_2, d_2,\cdots,d_n, b_{n+1}]$ be the $2$-bridge knot given in section 1, then the Kauffman polynomial of $D$ is
$$<D>= (0,0,\d)B(b_{n+1})D(d_n)\cdots B(b_{i+1})D(d_i)\cdots B(b_{2})D(d_1)B_1(b_1) \left(\begin{matrix}1\\ 1\\ 1 \end{matrix}\right )$$
where $B_1(b_1), B(b_{i}), D(d_i)$ are matrices defined in the
previous section as

\noindent $B_1(b_1)=(b_{ij}), 1\leq i,j\leq 3,\ 
with\  
b_{ij}=\left\{\begin{array}{cc}
m_{11}s^{b_1}+m_{12}(-s^{-1})^{b_1}+m_{13}\a^{-b_1}  & \mbox{if $i=j=1$}\\ & \\
m_{21}s^{b_1}+m_{22}(-s^{-1})^{b_1}+m_{23}\a^{-b_1} & \mbox{if $i=j=2$}\\ &\\
m_{31}s^{b_1}+m_{32}(-s^{-1})^{b_1}+m_{33}\a^{-b_1} & \mbox{if $i=j=3$} \\ &\\
0 & \mbox{Otherwise} \\
\end{array}
\right.$
\vspace{0.2cm}

\noindent $B(b_i)= M\left(\begin{matrix}s^{b_i} & 0 & 0\\ 0 & (-s^{-1})^{b_i} & 0\\ 0 & 0 & \a^{-b_i} \end{matrix}\right )M,$
and 
$D(d_i)= \left(\begin{matrix}s^{d_i} & 0 & 0\\ 0 & (-s^{-1})^{d_i} & 0\\ 0 & 0 & \a^{-d_i} \end{matrix}\right ).$

\end{theorem}

\begin{proof}
Using the linear maps defined in the previous section and their compositions, the $2$-bridge
knot
$D=C\circ B(b_{n+1})\circ D(d_n)\circ\cdots \circ B(b_{i+1})\circ D(d_i)\circ\cdots \circ B(b_{2})\circ D(d_1)\circ B_1(b_1)(\quad\raisebox{-3mm}{\epsfxsize0in\epsffile{parra2.ai}}\quad).$

As each of these maps is a linear transformation between vector
spaces, it can be represented by its matrix with respect to the
basis $\{{e_1}_h, {e_2}_h,{e_3}_h\}$. Note $\quad\raisebox{-3mm}{\epsfxsize0in\epsffile{parra2.ai}}\quad= {e_1}_h+
{e_2}_h +{e_3}_h$, so

$D=C\circ B(b_{n+1})\circ D(d_n)\circ\cdots \circ B(b_{i+1})\circ D(d_i)\circ\cdots \circ B(b_{2})\circ D(d_1)\circ B_1(b_1)({e_1}_h+
{e_2}_h +{e_3}_h),$

$=C\circ B(b_{n+1})\circ D(d_n)\circ\cdots \circ B(b_{i+1})\circ D(d_i)\circ\cdots \circ B(b_{2})\circ D(d_1) \left( \left(\begin{matrix}{e_1}_h & {e_2}_h & {e_2}_h \end{matrix}\right ) B_1(b_1)\left(\begin{matrix}1\\ 1\\ 1 \end{matrix}\right )\right)$, in matrices,

$=C\left( \left(\begin{matrix}{e_1}_h & {e_2}_h & {e_2}_h \end{matrix}\right ) B(b_{n+1}) D(d_n)\cdots  B(b_{i+1}) D(d_i)\cdots  B(b_{2}) D(d_1)B_1(b_1)\left(\begin{matrix}1\\ 1\\ 1 \end{matrix}\right )\right)$

$=(0,
0,\d\emptyset )B(b_{n+1}) D(d_n) \cdots   B(b_{i+1}) D(d_i) \cdots   B(b_{2})  D(d_1) B_1(b_1)
\left(\begin{matrix}1\\ 1\\ 1 \end{matrix}\right )$ \newline by Lemma 4.

Take the Kauffman polynomial, we have
$$<D>=(0,0,\d)B(b_{i+1})D(d_n)\cdots B(b_{i+1})D(d_i)\cdots B(b_{2})D(d_1)B_1(b_1) \left(\begin{matrix}1\\ 1\\ 1 \end{matrix}\right )$$
since the Kauffman polynomial of the empty link $\emptyset$ is $<\emptyset>=1$.
\end{proof}

\section{An Example--The Kauffman polynomial of the Whitehead link}

Here we demonstrate how to calculate the Kauffman polynomial of
the Whitehead link using linear maps and matrices. We choose the diagram $W$ with the
continued fraction notation $F(W)=[-2,1,-2]$ for the Whitehead link, then the corresponding matrices are as follows,
$$B_1(-2)=(b_{ij}), 1\leq i,j\leq 3,\ 
where\  
b_{ij}=\left\{\begin{array}{cc}
m_{11}s^{-2}+m_{12}s^{2}+m_{13}\a^{2}  & \mbox{if $i=j=1$}\\ & \\
m_{21}s^{-2}+m_{22}s^{2}+m_{23}\a^{2} & \mbox{if $i=j=2$}\\ &\\
m_{31}s^{-2}+m_{32}s^{2}+m_{33}\a^{2} & \mbox{if $i=j=3$} \\ &\\
0 & \mbox{Otherwise} \\
\end{array}
\right.$$

$D(1)=\left(\begin{matrix}s & 0 & 0\\ 0 & -s^{-1} & 0\\ 0 & 0 & \a^{-1} \end{matrix}\right ),$ and
$B(-2)=M\left(\begin{matrix}s^{-2} & 0 & 0\\ 0 & s^{2} & 0\\ 0 & 0 & \a^{2} \end{matrix}\right )M,$

According to Theorem 1, the Kauffman polynomial of $W$ is

$<W>= (0,0,\d)B(-2)D(1)B_{1}(-2)\left(\begin{matrix}1\\ 1\\ 1 \end{matrix}\right )$

$=\displaystyle{\frac{1}{\a^2\d s^4}}(\a^7s^4-2\a^4(-1 + \d^2)s^3(-1+s^2) + \a^3(1-s^2 + s^4)(1+(-1+2\d)s^2+s^4)-(-1+\d)\d s(-1+2s^2-2s^4+s^6) + \a^2(-2+\d+\d^2)s(-1+ 2s^2 -2s^4 + s^6) -2 \a^5(s^2 -s^4+ s^6) + \a (s^2(-1 + s^2)^2- \d^2 s^2 (-1 + s^2)^2 + \d^3 (s^2 -s^4+ s^6) + \d(-1 + s^2-2 s^4 + s^6- s^8 )))$.
Substitute in $\delta=\left({\dfrac{\alpha-\alpha^{-1}} {\ s - \
s^{-1}}}+1\right)$,
$<W>=\displaystyle{\frac{1}{\a^3s^4(-1+s^2)^2}}(-\a^2+ \a^4 - \a^3s + \a^5s + s^2 + \a^2 s^2 -\a^4s^2 - \a^6s^2 + \a s^3 + 2\a^3s^3 - 2\a^5s^3 - \a^7s^3 - 2s^4 + \a^4s^4 + 2\a^6s^4 - \a s^5 - 3\a^3s^5 + \a^5 s^5 + 3\a^7s^5 + 3s^6 - 2\a^2s^6 - \a^4s^6 - 2\a^6s^6 + \a s^7 + 3\a^3 s^7 - \a^5s^7 - 3\a^7s^7 - 2s^8 + \a^4 s^8 + 2\a^6s^8 -\a s^9 - 2\a^3 s^9 + 2\a^5 s^9 + \a^7s^9 + s^{10} + \a^2s^{10} - \a^4s^{10} - \a^6s^{10} + \a^3 s^{11} - \a^5 s^{11} - \a^2s^{12} + \a^4 s^{12}).$ 

Our calculations are carried out using Mathematica.

\end{document}